\documentclass[reqno]{amsart}
\usepackage{latexsym}
\usepackage{amssymb}
\newtheorem{Theo}{Theorem}
\newtheorem{Prop}[Theo]{Proposition}

\begin{document}
\title{The order of elements in Sylow $p$-subgroups of the symmetric
  group}
\author{Jan-Christoph Schlage-Puchta}
\begin{abstract}
Define a random variable $\xi_n$ by choosing a conjugacy class $C$ of
the Sylow $p$-subgroup of $S_{p^n}$ by random, and let $\xi_n$ be the
logarithm of the order of an element in $C$. We show that $\xi_n$ has
bounded variance and mean order $\frac{\log n}{\log p}+O(1)$, which differs
greatly from the average order of elements chosen with equal
probability.
\end{abstract}
\maketitle
In a sequence of papers (\cite{ET1}--\cite{ET7}), P. Erd{\H o}s and
P. Tur{\'a}n developed a statistical theory of the symmetric group $S_n$
on $n$ letters. Here a statistical theory describes properties of
almost all elements of a sequence of groups, or statements on
densities of certain elements or other associated structures,
e.g. characters.

P. Tur{\'a}n posed the problem of developing a statistical theory for
subgroups of $S_n$, in particular for Sylow subgroups. This was done
by P. P. P{\'a}lfy and M. Szalay (\cite{PS1}, \cite{PS2}, \cite{PS3}).

In particular, in the second of these papers they proved the
following:
\begin{Theo} Define the random variable $\zeta_n$ as follows: Choose an element
  $g$ of $P_n$, the Sylow $p$-subgroup of $S_{p^n}$, by random, and
  define $\zeta_n=\log_p o(g)$, where $o(g)$ is the order of $g$. Then
  $\zeta_n$ has bounded variance, and there are positive constants $c_1,
  c_2$, such that $c_1n<M_n<c_2n$, where $M_n$ denotes the mean value
  of $\zeta_n$.
\end{Theo}
Recently M. Ab{\'e}rt and B. Vir{\'a}g \cite{AV} showed, that the mean value
is in fact asymptotically equal to $c_pn$, where $c_p$ is given as
solution of the equation
\[
\big((1-c_p)/c_p\big)\log(1-c_p) + \log c_p = \log (1-1/p).
\]
In this note the analogous question is studied, however, we choose a
conjugacy class instead of an element.

Define $h_n(k)$ to be the number of conjugacy classes of $P_n$
consisting of elements of order $\leq p^k$. Our approach is based on the
following Theorem, the second half of which was proven by P. P. P{\'a}lfy
and M. Szalay \cite{PS1}.
\begin{Theo}
For $n, k\geq 1$ we have the recurrence relation
\begin{equation}
\label{eq:Rec}
h_{n+1}(k)=\frac{1}{p}\big(h_n(k)^p-h_n(k)\big)+h_n(k)+(p-1)h_n(k-1).
\end{equation}
In particular, for $h_n(n)$, the total number of conjugacy classes, we
have
\[
h_n(n)=\left[p^{\frac{\gamma+1}{p-1}}|P_n|^\gamma\right]-\delta_p,
\]
where $0<\gamma<1$ is a constant depending on $p$, and $\delta_p=1$ if $p=2$,
and $0$ otherwise.
\end{Theo}
Using this recurrence relation we will prove the following
theorem. Define $\xi_n$ to be $\log_po(C)$, where $C$ is a conjugacy
class chosen among the conjugacy classes of $P_n$ at random.
\begin{Theo}
The random variable $\xi_n$ has mean value $\frac{\log n}{\log p}+O(1)$
and has bounded variance.
\end{Theo}
Comparing Theorem 1 and Theorem 3, one sees that there are few
conjugacy classes of order $\sim c_pn$, containing almost all elements,
whereas almost all conjugacy classes are of order $\sim\frac{\log n}{\log
  p}$ but contain a neglectable proportion of all elements. It is
interesting to note that P. Erd{\H o}s and P. Tur{\'a}n discovered the
converse phenomenon in the symmetric group $S_n$: Almost all elements
of $S_n$ have order $e^{(1/2+o(1))\log^2 n}$, whereas almost all
conjugacy classes have order $e^{(c+o(1))\sqrt{n}}$, where $c$ is some
positive constant.

To prove Theorem 2, note that $P_{n+1}=P_{n}\wr C_p$, where $C_p$ is the
cyclic group of order $p$. To count the conjugacy classes in
$P_{n+1}$, we distinguish three cases: Conjugacy classes, which are
contained in the base group $P_n\times\cdots\times P_n$, and which are embedded in
the diagonal of this product, conjugacy classes which are
non-diagonally contained in the base group, and conjugacy classes,
which are not contained in the base group. Clearly, these cases are
disjoint and cover all possibilities. The conjugacy classes
of the first kind correspond to conjugacy classes of $P_n$, and the
action of $C_p$ leaves these conjugacy classes invariant, hence, there
are $h_n(k)$ conjugacy classes of elements of order $\leq p^k$ of the first
kind. There are $h_n(k)^p-h_n(k)$ conjugacy classes in $P_n\times\dots\times
P_n$ which do not contain diagonal elements, and the action of
$C_p$ does not leave any such class invariant, hence, there are
$\frac{1}{p}(h_n(k)^p-h_n(k))$ conjugacy classes of the second
kind. To count conjugacy classes of the third kind, write an element
of $P_{n+1}$ as $(f, \sigma)$, where $f:\{1, \ldots, p\}\to P_n$, and
$\sigma\in C_p$. Conjugation cannot change $\sigma$, hence, the total number of
conjugacy classes of the third kind equals $(p-1)$ times the number
of conjugacy classes of elements having $\sigma(1)=2$. For an element $x\in
P_n$ and an integer $i\leq p$, define the function $g$ by $g(i)=x$,
$g(j)=1$ for $j\neq i$. Then we find $(f, \sigma)^{(g, \mathrm{id})} = (h,
\sigma)$, where
\[
h(j) = \begin{cases} f(j), & j\neq i, i+1\\
 xf(i), & j=i\\ f(i+1)x^{-1}, & j=i+1\end{cases}.
\]
In particular, every element of the third kind is conjugated to an
element $(f, \sigma)$ with $f(i)=1$ for all $i\neq 1$. Moreover, choosing
$g(i)=x$ for all $i\in\{1, \ldots, p\}$, we find that conjugacy classes of the
third kind correspond to conjugacy classes of $P_n$. Finally, the
order of $(f, \sigma)$ in $P_{n+1}$ equals $p$ times the order of
$\prod_{i=1}^p f(i)$ in $P_n$, hence, there are $h_n(k-1)$ conjugacy
classes of elements of order $k$ in $P_{n+1}$ of the third
kind. Adding up the contribution of these three classes yields the
recurrence relation (\ref{eq:Rec}). The formula for $h_n(n)$ follows
from this by induction. A different proof of this formula counting
characters was given by P. P. P{\'a}lfy and M. Szalay \cite{PS1}.

We begin the proof of Theorem 3 with two simple remarks: First it is
obvious from 
the definition that for fixed $n$, $h_n(k)$ is increasing with
$k$. Further note that $h_{n+1}(1)>h_n(1)+(p-1)$, thus, for $n\nearrow\infty$,
$h_n(1)$ tends to $\infty$.

We will write
\[
\alpha_n(k) := \frac{h_n(k-1)}{h_n(n)}
\]
and
\[
\beta_n(k) := \frac{h_n(k-1)}{h_n(k)}
\left(\frac{h_{n-1}(k)}{h_{n-1}(k-1)}\right)^p.
\]
We first collect some properties of $\alpha$ and $\beta$.
\begin{Prop}
Define $\alpha$ and $\beta$ as above, and assume that $1\leq k\leq n$.
\begin{enumerate}
\item We have
\begin{equation}
\label{eq:Prop1}
\beta_n(k) =
\frac{1+(p-1)h_{n-1}(k-1)^{1-p}+p(p-1)h_{n-1}(k-2)h_{n-1}(k-1)^{-p}}
{1+(p-1)h_{n-1}(k)^{1-p}+p(p-1)h_{n-1}(k-1)h_{n-1}(k)^{-p}}
\end{equation}
\item We have $|\log \beta_n(k)|\ll h_{n-1}(k-1)^{1-p}$.
\item We have
\[
\alpha_n(k) = \beta_k(k)^{p^{n-k}} (\beta_{k+1}(k)\beta_{k+1}(k+1))^{p^{n-k-1}}\cdots
(\beta_n(k)\cdots\beta_n(n)).
\]
\item We have $\beta_n(n)=1-\frac{(p-1)^n}{h_n(n)}$. Furthermore, for
  $p=2$, we have
\[
\beta_n(n-1) = 1+O\left(\frac{h_{n-2}(n-2)}{h_n(n)}\right).
\]
\item As $k$ tends to $\infty$, the asymptotic $\log \alpha_n(k)\sim
  -p^{-n-k}\frac{(p-1)^k}{h_k(k)}$ holds.
\end{enumerate}
\end{Prop}

Before we prove this, we indicate how Theorem 3 follows from
Proposition 4. Obviously, for fixed $n$, $\alpha_n(k)$ is increasing with
$k$, $\alpha_n(1)=\frac{1}{h_n(n)}\to 0$, and $\alpha_n(n+1)=1$, hence, for $n\geq 1$
  there is a unique integer $k_0(n)$ with the property
  $\alpha_n(k_0)\leq1/2<\alpha_n(k_o+1)$. By 4.5 we get
\[
-\log 2\geq\log\alpha_n(k_0)\sim -p^{n-k_0}\frac{(p-1)^{k_0}}{h_{k_0}(k_0)}
\]
Taking logarithms again and using the value for $h_{k_0}(k_0)$ given
by Theorem 1, we get
\[
\log\log 2\leq (1+o(1))\Big((n-k_0)\log p + k_0\log(p-1)
-\frac{\gamma+1}{p-1}\log p-\gamma\frac{p^{k_0}-1}{p-1}\log p\Big).
\]
From this inequality it is easy to see that $k_0\leq\frac{\log n}{\log
p}+O(1)$. In the same way, starting from the inequality
$\alpha_n(k_0+1)>1/2$, we get $k_0\geq\frac{\log n}{\log p}+O(1)$, thus
$k_0=\frac{\log n}{\log p}+O(1)$. We want to show that $k_0$ is
close to the mean value of $\xi_n$, and that the variance of $\xi_n$ is
bounded. Both statements follow at once, if we can show that the
mean value of $(\xi_n-k_0)^2$ is bounded. We estimate this
value as follows.
\begin{eqnarray*}
\mathbf{E}(\xi_n-k_0)^2 & = &
\frac{1}{h_n(n)}\sum_{k=0}^n(h_n(k)-h_n(k-1))(k-k_0)^2\\
& = & \sum_{k=0}^n(\alpha_n(k+1)-\alpha_n(k))(k-k_0)^2\\
& \leq & \sum_{k=0}^{k_0}\alpha_n(k+1)(k-k_0)^2 +
\sum_{k=k_0+1}^{n}(1-\alpha_n(k))(k-k_0)^2.
\end{eqnarray*}
To estimate these sums, we note first that $\alpha_n(k)^2>\alpha_n(k-1)$ for $k$
greater than some absolute constant, for by 4.5 we have
\[
\frac{\log\alpha_n(k-1)}{\log\alpha_n(k)} =
(1+o(1))\frac{p}{p-1}\frac{h_k(k)}{h_{k-1}(k-1)},
\]
and since the right hand side tends to $\infty$, it will eventually become
$\geq 2$. Using this fact together with the definition of $k_0$, we obtain the
estimates
\begin{eqnarray*}
\alpha_n(k_0-d) & < & 2^{-2^d},\\
\alpha_n(k_0+d) & > & \sqrt[2^{d-1}]{1/2} = 1-\frac{\log 2}{2^{d-1}} +
O(2^{-2d}).
\end{eqnarray*}
Hence, both sums in the estimate for $\mathbf{E}(\xi_n-k_0)^2$ can be
estimated by converging sums of the form $d^22^{-d}$, and we obtain
$\mathbf{E}(\xi_n-k_0)^2\ll 1$, which proves Theorem 3.

Therefore, it remains to prove Proposition 4.

1. In the definition of $\beta_n(k)$, replace the values of $h$ in the
first fraction with the recurrence relation (\ref{eq:Rec}) to obtain
\begin{eqnarray*}
\beta_n(k) & = &
\frac{h_{n-1}(k-1)^p+(p-1)h_{n-1}(k-1)+p(p-1)h_{n-1}(k-2)}
{h_{n-1}(k)^p+(p-1)h_{n-1}(k)+p(p-1)h_{n-1}(k-1)}
\left(\frac{h_{n-1}(k)}{h_{n-1}(k-1)}\right)^p\\
 & = &
 \frac{1+(p-1)h_{n-1}(k-1)^{1-p}+p(p-1)h_{n-1}(k-2)h_{n-1}(k-1)^{-p}}
{1+(p-1)h_{n-1}(k)^{1-p}+p(p-1)h_{n-1}(k-1)h_{n-1}(k)^{-p}}.
\end{eqnarray*}
2. We use the (\ref{eq:Prop1}). To give an upper bound, we estimate the
denominator by 1. In the numerator we replace $h_{n-1}(k-2)$ by
$h_{n-1}(k-1)$, which increases the fraction, too. Thus we obtain
the upper bound $1+(p^2-1)h_{n-1}(k-1)^{1-p}$. In the same way we
obtain a lower bound, and both bounds together yield a bound for
$|\log\beta_n(k)|$. 

3. This statement follows from the definition of $\alpha$ by a simple
computation:
\begin{eqnarray*}
\frac{h_n(k-1)}{h_n(k)} & = &
\beta_n(k)\left(\frac{h_{n-1}(k-1)}{h_{n-1}(k)}\right)^p\\
 & = & \beta_n(k)\beta_{n-1}(k)^p\cdots \beta_k(k)^{p^{n-k}}\\
\alpha_n(k) & = &
\frac{h_n(k-1)}{h_n(k)}\frac{h_n(k)}{h_n(k+1)}\cdots\frac{h_n(n-1)}{h_n(n)}\\
 & = & \big(\beta_n(k)\beta_{n-1}(k)^p\cdots \beta_k(k)^{p^{n-k}}\big)\\
&&\quad \big(\beta_n(k+1)\beta_{n-1}(k+1)^p\cdots \beta_{k+1}(k+1)^{p^{n-k-1}}\big)\\
&&\quad\cdots\, \big(\beta_n(n-1)\beta_{n-1}(n-1)^p\big)\beta_n(n).
\end{eqnarray*}
Now rearranging terms according to the exponent proves our claim.

4. To compute $\beta_n(n)$, it suffices to compute $h_n(n-1)$. We have
\begin{eqnarray*}
h_n(n-1) & = &
\frac{1}{p}\big(h_{n-1}(n-1)^p-h_{n-1}(n-1)\big)+h_{n-1}(n-1)+(p-1)h_{n-1}(n-2)\\
& = & h_n(n) - (p-1)\big(h_{n-1}(n-1)-h_{n-1}(n-2)\big). 
\end{eqnarray*}
From this we deduce by induction that
$h_n(n)-h_n(n-1)=(p-1)^{n-1}\big(h_1(1)-h_1(0)\big)$. But
$h_1(1)-h_1(0)=p-1$, since $P_1$ is cyclic of order $p$, and we get
\[
\beta_n(n)=\frac{h_n(n-1)}{h_n(n)}=1-\frac{(p-1)^n}{h_n(n)}.
\]
Now assume that $p=2$. Then we have
\begin{eqnarray*}
h_n(n-2) & = &
\frac{1}{2}\big(h_{n-1}(n-2)^2-h_{n-1}(n-2)\big)+h_{n-1}(n-2)+h_{n-1}(n-3)\\
 & = &
 \frac{1}{2}\Big(\big(h_{n-1}(n-1)-1\big)^2-\big(h_{n-1}(n-1)-1\big)\Big)+
h_{n-1}(n-1)-1+h_{n-1}(n-3)\\
 & = & h_n(n-1)-h_{n-1}(n-1)-h_{n-1}(n-2)+h_{n-1}(n-3)\\
 & = & h_n(n)-2h_{n-1}(n-1)+h_{n-1}(n-3),
\end{eqnarray*}
which implies $h_n(n)-h_n(n-2)\ll h_{n-1}(n-1)$; using this estimate to
bound $h_{n-1}(n-1)-h_{n-1}(n-3)$, we obtain
From this we obtain
\[
h_n(n-2)= h_n(n)-h_{n-1}(n-1)+O(h_{n-2}(n-2)).
\]
Now we can compute $\beta_n(n-1)$:
\begin{eqnarray*}
\beta_n(n-1) & = &
\frac{h_n(n-2)}{h_n(n-1)}\left(\frac{h_{n-1}(n-1)}{h_{n-1}(n-2)}\right)^2\\
 & = & \frac{h_n(n)-h_{n-1}(n-1)+O(h_{n-2}(n-2))}{h_n(n)-1}\cdot
\left(\frac{h_{n-1}(n-1)}{h_{n-1}(n-1)-1}\right)^2\\
 & = & 1-\frac{h_{n-1}(n-1)}{h_n(n)} +
 O\left(\frac{h_{n-2}(n-2)}{h_n(n)}\right) + \frac{2}{h_{n-1}(n-1)} +
 O\left(\frac{1}{h_{n-1}(n-1)^2}\right)\\
 & = & 1+\frac{2h_n(n)-h_{n-1}(n-1)^2}{h_{n-1}(n-1)h_n(n)} +
 O\left(\frac{h_{n-2}(n-2)}{h_n(n)}\right)\\
 & = & 1+O\left(\frac{h_{n-2}(n-2)}{h_n(n)}\right),
\end{eqnarray*}
where in the last line we used the recurrence relation in the form
\[
h_n(n)=\frac{1}{2}h_{n-1}(n-1)^2+O(h_{n-1}(n-1)).
\]

5. By 4.3 we can express $\alpha$ in terms of $\beta$, and then we will use 4.2
and 4.4 to estimate the resulting expression. We have
\begin{eqnarray}
\log\alpha_n(k) & = & p^{n-k}\log\beta_k(k) +
p^{n-k-1}(\log\beta_{k+1}(k)+\log\beta_{k+1}(k+1)) +\nonumber\\
&&\quad\dots+(\log\beta_{n}(k)+\ldots+\log\beta_{n}(n))\label{eq:Prop5}\\
 & = & -p^{n-k}\frac{(p-1)^k}{h_k(k)} + O\big(p^{n-k}(p-1)^{2k}k_k(k)^{-2}\big)+
O\left(\sum_{\kappa=k+1}^np^{n-\kappa}\sum_{\nu=k}^\kappa\log\beta_\kappa(\nu)\right). \nonumber
\end{eqnarray}
The first error term is of lesser order than the main term, provided
that $k\nearrow\infty$. To bound the second error term, we first consider the sum
over the range $k+2\leq\kappa\leq n$. Using 4.2, we find
\begin{eqnarray*}
\sum_{\kappa=k+2}^np^{n-\kappa}\sum_{\nu=k}^\kappa\log\beta_\kappa(\nu) & \ll & \sum_{\kappa=k+2}^np^{n-\kappa}\sum_{\nu=k}^\kappa
h_{\kappa-1}(\nu-1)^{1-p}\\
 &\ll & \sum_{\kappa=k+2}^np^{n-\kappa} \frac{\kappa-k+1}{h_{\kappa-1}(k-1)}\\
 & \ll & \frac{p^{n-k}}{h_{k+1}(k-1)}.
\end{eqnarray*}
As can be seen from 4.2 and 4.3, we have $h_n(n-2)\sim h_n(n)$ as $n\to\infty$,
hence, the last quantity is of the
same order of magnitude as $\frac{p^{n-k}}{h_{k+1}(k+1)}$, which is negligible
compared to the main term provided that $k\nearrow\infty$. To bound the term
coming from $\kappa=k+1$, we 
distinguish the cases $p=2$ and $p\geq 3$. In the latter case, 4.2
implies that
\[
\log\beta_{k+1}(k)+\log\beta_{k+1}(k+1) \ll h_k(k-1)^{1-p}\leq h_k(k-1)^{-2}.
\]
As in the proof of 4.4 we find that
$h_k(k)-h_k(k-1)=(p-1)^k=o(h_k(k))$, hence, for $k\nearrow\infty$, 
$h_k(k-1)^{-2}=o(h_k(k)^{-1})$, and the contribution of this term is
negligible. If on the other hand $p=2$, we use 4.4 to deduce
\[
\log\beta_{k+1}(k)+\log\beta_{k+1}(k+1)\sim-\frac{1}{h_{k+1}(k+1)} +
O\left(\frac{h_{k-1}(k-1)}{h_{k+1}(k+1)}\right).
\]
Since $h_{k+1}(k+1)\gg h_k(k)^2\gg h_{k-1}(k-1)^4$, we obtain
\[
\log\beta_{k+1}(k)+\log\beta_{k+1}(k+1)\ll
\frac{1}{h_{k}(k)^2}+\frac{1}{h_{k}(k)^{3/2}} = o(h_k(k)^{-1}).
\]
Hence, the main term in (\ref{eq:Prop5}) dominates the error terms as
$k\nearrow\infty$, and 4.5 is proven.

\end{document}